\documentclass[10pt]{amsart}

\usepackage{amsmath,amsfonts,amssymb}
\usepackage{amsthm}
\newtheorem{theorem}{Theorem}[section]
\newtheorem{corollary}[theorem]{Corollary}
\newtheorem{proposition}[theorem]{Proposition}
\newtheorem{lemma}[theorem]{Lemma}
\newtheorem{remark}[theorem]{Remark}

\newcommand{\Ac}{\mathcal A}

\newcommand{\Vc}{\mathcal V}

\newcommand{\cM}{\mathcal M}
\newcommand{\Sc}{\mathcal S}
\newcommand{\Zc}{\mathcal Z}
\newcommand{\Nc}{\mathcal N}

\newcommand{\Hc}{\mathcal H}

\newcommand{\Lc}{\mathcal L}
\newcommand{\Jc}{\mathcal J}

\newcommand{\Rb}{\mathbb{R}}
\newcommand{\Cb}{\mathbb{C}}
\newcommand{\Hb}{\mathbb{H}} 
\newcommand{\Nb}{\mathbb{N}}

\newcommand{\Ad }{\mbox{\rm Ad}}

\newcommand{\Coad }{\mbox{\rm Coad}}
\newcommand{\Ind}{\mbox{\rm Ind}}

\newcommand{\Id}{\mbox{\rm Id}}

\newcommand{\nn}[1]{| #1 |}
\newcommand{\nd}[1]{\left\|#1 \right\|}

\newcommand{\proj}[2]{\mbox{pr}_{#1}(#2)}

\newcommand{\fbs}[1]{\Jc_{#1}}
\newcommand{\fb}[2]{\fbs{#1}(#2)}

\newcommand{\fls}[2]{\Lc_{#1,#2}}
\newcommand{\fl}[3]{\fls{#1}{#2}(#3)}

\newcommand{\flbs}[2]{\bar{\Lc}_{#1,#2}}
\newcommand{\flb}[3]{\flbs{#1}{#2}(#3)}

\begin{document}

\title[]
{The Bounded Spherical Functions for the Free  Two-step Nilpotent Lie Group}       

\author[]
{V\'eronique Fischer}     

\address{Department of Mathematics of G\"oteborg University,  S-412 96  G\"oteborg Sweden}
\email{veronfi@math.chalmers.se, v.fischer@sns.it}


\begin{abstract}
  In this paper, we give the expressions 
  for the bounded spherical functions, 
  or equivalently the spherical functions of positive type,
  for the free two-step nilpotent Lie groups 
  endowed  with the actions
  of orthogonal groups or their special subgroups.
  Next we deduce some results about the (Kohn) sub-Laplacian,
  and we compute the radial Plancherel measure.
\end{abstract}

\maketitle


\section{Introduction}
A (connected, simply connected)
nilpotent Lie group which forms with a compact Lie group 
a Gelfand pair, is at most of step two, 
and the bounded spherical functions are 
the spherical functions of positive type \cite{pgelf}.
The cases of the Heisenberg group with some subgroups
of the unitary matrix group  are well known \cite{bjr},
and the bounded spherical functions are then explicit.
In this paper, we are interested by  the Gelfand pair formed by 
the free two-step nilpotent Lie groups 
with the actions of orthogonal groups.
The expressions of some of the corresponding bounded spherical functions
were given
in \cite[Section 6]{stric}
with a sketched proof.
Here we give the expression of all such functions with a complete proof,
and we obtain the corresponding eigenvalues for the sub-Laplacian 
and the radial Plancherel measure. 

This paper is organized as follows.
After recalling some definitions and setting some notations, 
we give in the second section the statement of the main results:
the expressions of the bounded spherical functions
on the free two-step nilpotent Lie groups.
In the third section, 
we recall a few facts 
about spherical functions 
and representations,
which allow to construct  
our bounded spherical functions in the following section.
We also give an equivalent method of construction
from which we obtain some properties of the sub-Laplacian
and the radial Plancherel formula in the fifth section.\\
We shall omit some computations 
and the proof for the case of the special orthogonal group.
We refer the interested reader to the French thesis of the author
\cite{moi}. 


\section{The Free Two-step Nilpotent Lie Groups}

Here we give definitions and notations 
for the free two-step nilpotent Lie groups and algebras;
we also present the action of orthogonal groups.

\paragraph{First Definition.}
Let $\Nc_p$ be the (unique up to isomorphism)
free two-step nilpotent Lie algebra with $p$ generators.
The definition using the universal property of the free nilpotent Lie algebra can
be found in \cite[Chapter V \S 4]{jacobson}.
Roughly speaking, $\Nc_p$ is a (nilpotent) Lie algebra
with $p$ generators $X_1,\ldots,X_p$,
such that the vectors  $X_1,\ldots,X_p$ and
$X_{i,j}=[X_i,X_j], i<j$ form a basis;
we call this basis the canonical basis of $\Nc_p$. \\
We denote by
$\Vc$ and $\Zc$,
the vector spaces generated by 
the families of vectors $X_1,\ldots,X_p$
and $X_{i,j}:=[X_i,X_j], 1\leq i< j\leq p$
respectively;
these families become the canonical base of $\Vc$ and $\Zc$.
Thus $\Nc_p=\Vc\oplus\Zc$,
and $\Zc$ is the center of $\Nc_p$. 
With the canonical basis,
the vector space $\Zc$ can be identified with the vector space of
antisymmetric $p\times p$-matrices $\Ac_p$.
Let $z=\dim \Zc=p(p-1)/2$.\\
The connected simply connected nilpotent Lie group
which corresponds to $\Nc_p$ is called the
free two-step nilpotent Lie group and is denoted $N_p$.
We denote by $\exp :\Nc_p\rightarrow N_p$ the exponential map.\\
In the following, we use the notations $X+A\in\Nc$, $\exp(X+A)\in N$
when $X\in \Vc,A\in\Zc$.
We write $p=2p'$ or $2p'+1$.

\paragraph{A Realization of \mathversion{bold}{$\Nc_p$}.}
We now present here a realization of $\Nc_p$,
which will be helpful to define more naturally
the action of the orthogonal group and representations of $N_p$.\\
Let $(\Vc,<,>)$ be an Euclidean space with dimension $p$.
Let $O(\Vc)$ be the group of orthogonal transformations of $\Vc$,
and $SO(\Vc)$ its special subgroup.
Their common Lie algebra denoted by $\Zc$,
is identified with the vector space 
of antisymmetric transformations of $\Vc$.
Let $\Nc=\Vc\oplus\Zc$ 
be the exterior direct sum of the vector spaces $\Vc$ and $\Zc$.\\
Let $[,]:\Vc\times\Vc\rightarrow \Zc$ 
be the bilinear application given by :
$$
[X,Y].(V)\,=\, <X,V>Y-<Y,V>X
\quad X,Y,V\in \Vc
\quad.
$$
We also denote by $[,]$ the bilinear application
extended to $\Nc\times\Nc\rightarrow \Nc$ by:
$$
[.,.]_{\Nc\times \Zc}\,=\,
[.,.]_{\Zc\times \Nc}\,=\,
0\quad.
$$
This application is a Lie bracket.
It endows $\Nc$ with the structure 
of a two-step nilpotent Lie algebra.\\
As the elements $[X,Y]$, $X,Y\in \Vc$ generate the vector space
$\Zc$, we also define a scalar product $<,>$ on $\Zc$ by:
$$
<[X,Y],[X',Y']>
\,=\, 
<X,X'><Y,Y'>-<X,Y'><X',Y>,
\quad,
$$
where $X,Y,X',Y'\in \Vc$.\\
It is easy to see  $\Nc$ as a realization of $\Nc_p$
when an orthonormal basis
$X_1,\ldots,X_p$  of $(\Vc,<,>)$
is fixed.\\
We remark that
$<[X,Y],[X',Y']>= <[X,Y]X',Y'>$,
and so we have  for an antisymmetric transformation $A\in \Zc$,
and for $X,Y\in \Vc$:
\begin{equation}
  \label{equality_bracket_scpdt}
  <A,[X,Y]>
  \,=\,
  <A.X,Y>
  \quad.  
\end{equation}
This equality can also be proved directly using the canonical basis of $\Nc_p$.

\paragraph{Actions of Orthogonal Groups.}
We denote by $O(\Vc)$ the group of orthogonal linear maps of $(\Vc, <,>)$,
and by  $O_p$ the group of orthogonal $p\times p$-matrices.

\paragraph{On \mathversion{bold}{$\Nc_p$} and \mathversion{bold}{$N_p$}.}
The group $O(\Vc)$ acts 
on the one hand by automorphisms on $\Vc$,
on the other hand by the adjoint representation $\Ad_\Zc$ on $\Zc$.
We obtain an action of $O(\Vc)$ on $\Nc=\Vc\oplus\Zc$.
Let us prove that this action respects the Lie bracket of $\Nc$.
It suffices to show 
for $X,Y,V\in\Vc$ and $k\in O(\Vc)$ :
\begin{eqnarray*}
  [k.X,k.Y](V)
  &=&
  <k.X,V>k.Y-<k.Y,V>k.X\\
  &=&
  k.\left( <X, ^tk.V>Y-<Y, ^tk.V>X\right)\\
  &=&
  k.[X,Y](k^{-1}.V)
  \,=\, Ad_\Zc k.[X,Y]
  \quad.
\end{eqnarray*}
We then obtain that the group $O(\Vc)$ 
(and also its special subgroup $SO(\Vc)$)
acts by automorphism on the Lie algebra $\Nc$, 
and finally on the Lie group $N$.\\
Suppose an orthonormal basis  $X_1,\ldots,X_p$  of $(\Vc,<,>)$
is fixed;
then the vectors $X_{i,j}:=[X_i,X_j], 1\leq i< j\leq p$, 
form an orthonormal basis of $\Zc$ 
and we can identify:
\begin{itemize}
\item the vector space $\Zc$ with $\Ac_p$,
\item the group $O(\Vc)$ with $O_p$,
\item the adjoint representation $\Ad_\Zc$
  with the conjugate action of $O_p$ on $\Ac_p$:
  $k.A=kAk^{-1}$, where $k\in O_p$, $A\in \Ac_p$.
\end{itemize}
Thus the group $O_p\sim O(\Vc)$ acts on $\Vc\sim\Rb^p$ and $\Zc\sim\Ac_p$,
and consequently on $\Nc_p$.
Those actions can be directly defined;
and the equality $[k.X,k.Y]=k.[X,Y]$,
$k\in O_p$, $X,Y\in\Vc$, can then be computed.

\paragraph{On \mathversion{bold}{$\Ac_p$}.}
Now we describe the orbits of the conjugate actions 
of $O_p$ and $SO_p$ on $\Ac_p$.\\
An arbitrary antisymmetric matrix $A\in \Ac_p$ is $O_p$-conjugated 
to an antisymmetric matrix $D_2(\Lambda)$ where $      \Lambda=(\lambda_1,\ldots,\lambda_{p'})\in\Rb^{p'}$ and:
\begin{equation}
  \label{def_D2Lambda}
  D_2(\Lambda)
  \,:=\,
  \left[
    \begin{array}{cccc}
      \lambda_1 J &&\\
      0 & \ddots & 0&\\
      &&\lambda_{p'} J&\\
      &&&(0)
    \end{array}\right]
  \qquad \mbox{where}\quad
  J:=
  \left[  \begin{array}{cc}
      0&1\\
      -1&0
    \end{array}\right]
  \quad.
\end{equation}
($(0)$ means that a zero appears only in the case $p=2p'+1$.)
Furthermore, we can assume that $\Lambda$ is in $\bar{\Lc}$,
where we denote by $\bar{\Lc}$ the set of 
$\Lambda=(\lambda_1,\ldots,\lambda_{p'})\in\Rb^{p'}$,
such that
$\lambda_1\,\geq\, \ldots \,\geq\, \lambda_{p'}\,\geq\,0$.\\
An arbitrary antisymmetric matrix $A\in \Ac_p$ is $SO_p$-conjugated 
to 
$$
D_2^\epsilon(\Lambda)
:=
D_2(\lambda_1,\ldots,\lambda_{p'-1},\epsilon\lambda_{p'})
\quad,
$$
where $\epsilon=\pm1$
and $\Lambda=(\lambda_1,\ldots,\lambda_{p'})\in\bar{\Lc}$.


\section{Notations and Main Results}

We give here the notations for special functions
that will be used to present the main results of this paper.
First, we recall the definitions of the Bessel and Laguerre functions,
and we set some notations for parameters.
We give then the expression of the bounded spherical functions.

\paragraph{Notations for special functions.}
We will use the following well known functions :
\begin{itemize}
\item the Gamma function $\Gamma$,
\item the Laguerre polynomial 
  of type $\alpha$ and degree $n$:  $L_n^{\alpha}$ \cite[\S 5.1]{szego},
\item the Bessel function of type $\alpha$:  $J_\alpha$ 
  \cite[\S 1.71]{szego}, \cite[ch. II, I.1]{faraut_h}.
\end{itemize}
Let us now define the normalized Laguerre function  $\flbs n \alpha = \fls  {n}{\alpha} / C_{n+\alpha}^n $ where:
$$
\fl{n}{\alpha}{x}:=L_n^{\alpha}(x)e^{-\frac{x}{2}}
\quad\mbox{and}\quad 
C_{n+\alpha}^n=\frac{\Gamma(n+\alpha+1)}{n!\Gamma(\alpha+1)}
\quad,
$$
and the reduced Bessel function $\fbs \alpha$ by:
$$
\fb \alpha{z}:=
\Gamma (\alpha+1){(z/2)}^{-\alpha}J_\alpha(z)
\quad.
$$
Let $n = 1,2,\ldots$,
and let $dk$ denote the Haar probability measure
of the compact group $K=O_n$ or $SO_n$.
If $<,>$ denotes the Euclidean scalar product, 
and $\nn{.}$ the Euclidean norm of $\Rb^n$,
we recall for any fixed $x_0\in\Rb^n$ such that $\nn{x_0}=1$
\cite[ch.II, I.1]{faraut_h}:
\begin{equation}
  \label{formula_bessel}
  \fb {\frac {n-2}2}{\nn{x}}
  \,=\,
  \int_{K} e^{i<k.x,x_0>} dk
  \quad.  
\end{equation}

\paragraph{Parameters.}
To each  $\Lambda\in \bar{\Lc}$,
we associate:
$p_0$ the number of $\lambda_i\not=0$,
$p_1$ the number of distinct $\lambda_i\not=0$,
and $\mu_1,\ldots,\mu_{p_1}$ such that:
\begin{equation}
  \label{formula_mu_lambda}
  \{\mu_1>\mu_2>\ldots>\mu_{p_1}>0\}=
  \{\lambda_1\geq\lambda_2\geq\ldots\geq \lambda_{p_0}>0\} \quad.
\end{equation}
We denote by $m_j$ the number of $\lambda_i$ 
such that $\lambda_i=\mu_j$,
and we put:
\begin{equation}
  \label{formula_m_m'}
  m_0:=m'_0:=0
  \quad\mbox{and for}\;
  j=1,\ldots p_1\,\quad
  m'_j:=m_1+\ldots+m_j\quad.
\end{equation}
For $j=1,\ldots p_1$,
let $\mbox{pr}_j$ be the orthogonal projection of $\Vc$ onto the space 
generated by the vectors
$X_{2i-1},X_{2i}$, for
$i=m'_{j-1}+1,\ldots,m'_j$.\\
Let $\cM$ be the set of $(r,\Lambda)$
where $\Lambda\in \Lc$,
and $r\geq 0$,  
such that $r=0$ if $2p_0=p$.

\paragraph{Expression of the bounded spherical functions.}
The bounded spherical functions of $(N_p , K)$
for $K=O_p$ or $SO_p$,
are parameterized by 
\begin{itemize}
\item $(r,\Lambda)\in \cM$ 
  (with the previous notations $p_0$, $p_1$, $\mu_i$, $\mbox{pr}_j$ associated to $\Lambda$),
\item $l\in\Nb^{p_1}$ if $\Lambda\not=0$, otherwise $\emptyset$,
\item $\epsilon=\pm1$ if $K=SO_p$, otherwise  $\emptyset$.
\end{itemize}
Let $(r,\Lambda)$, 
$l$ and $\epsilon$
be such parameters.
If $\Lambda\not=0$,
we define  
the function
$\phi^{r,\Lambda,l,\epsilon}$ by:
\begin{equation}
  \label{def_sph_fcn_lag}
  \phi^{r,\Lambda,l,\epsilon}(n)
  \,=\,
  \int_{K}
  \Theta^{r,\Lambda,l,\epsilon}(k.n)
  dk,
  \quad n\in N_p
  \quad,  
\end{equation}
where $\Theta^{r,\Lambda,l,\epsilon}$ is given by:
\begin{equation}
  \label{def_theta_fcn}
  \Theta^{r,\Lambda,l,\epsilon}(\exp (X+A))
  \,=\,
  e^{i <rX_p,X>}
  e^{i<D_2^\epsilon(\Lambda),A>}
  \overset{p_1}
  {\underset{j=1}\Pi}
  \flb {l_j} {m_j-1} {\frac{\mu_j}2 \nn{\proj {j} X }^2}
  \quad.  
\end{equation}
If $\Lambda=0$, we define the function $\phi^{r,0}$ by:
\begin{equation}
  \label{def_sph_fcn_bes}
  \phi^{r,0}(n)
  \,=\,
  \fb {\frac{p-2}2}{r\nn{X}}
  \quad ,
  \quad n=\exp (X+A)\in N_p
  \quad.
\end{equation}
In Section~\ref{sec_proof},
we shall prove the following result in the case $K=O_p$
(the case $K=SO_p$ is similar and can be found in \cite{moi}):
\begin{theorem}
  \label{mainthm}
  The bounded spherical functions of $(N_p , K)$,
  for $K=SO_p$ or $K=O_p$,
  are the functions
  $\phi^{r,\Lambda,l,\epsilon}$
  given by~(\ref{def_sph_fcn_lag}) and~(\ref{def_sph_fcn_bes}),
  where $(r,\Lambda)\in \cM$ 
  and $l\in\Nb^{p_1}$ if $\Lambda\not=0$, 
  and $\epsilon=\pm1$ if $K=SO_p$,
  $\epsilon=\emptyset$ if $K=O_p$.
\end{theorem}
In Section~\ref{sec_rq_repN},
we shall also express $\phi^{r,0}$
and $\phi^{r,\Lambda,l}$
in terms of representations of $N_p$
and obtain their eigenvalues for the sub-Laplacian
and the radial Plancherel measure.


\section{Spherical Function and Representation}

In this section, we recall some of the properties of spherical functions,
Gelfand pairs and representations,
which will be used in the proof of Theorem~\ref{mainthm}.

In this article, we use the following conventions.\\
The semi-direct product $K\ltimes N$
of two groups $K$ and $N$ such that $K$ acts on $N$ by automorphism,
is defined by the law:
$$
(k_1,n_1),(k_2,n_2)\in K\ltimes N,
\quad
(k_1,n_1).(k_2,n_2)
\,=\,
(k_1k_2,n_1\, k_1.n_2)
\quad .
$$
All the groups are supposed locally compact, second countable and separable,
and their continuous unitary representations 
on separable Hilbert spaces.\\
For such a group $G$,
we denote by  $\hat{G}$ the quotient set 
of the irreducible representations
by the equivalence relation $\sim$.
We often identify a representation with its equivalence class.

\paragraph{Definitions and properties.}
Let $K$ be a compact group, which acts continuously on a group $N$.
Let $dn$ be a Haar measure on $N$,
and $dk$ the normalized Haar measure on $K$.
We assume that $dndk$ is a Haar measure on the group $G=K\ltimes N$, 
and that this group is unimodular.\\
Let $C^\natural (N)$ be the set 
of continuous compactly supported $K$-invariant functions on $N$.\\
A $K$-invariant function $\phi$ on $N$
is spherical on $N$ if for all $ f,g\in C^\natural(N)$ we have:
$$
\int_N f*g(n) \phi(n) dn
\,=\,
\int_N f(n) \phi(n) dn
\int_N g(n) \phi(n) dn
\;.
$$
If $\phi$ is a spherical function on $N$ for $K$,
the function $\Phi$ on $G=K\ltimes N$ given by
$\Phi(k,n)=\phi(n)$ is also called a spherical function 
on $G$ for $K$. 
\begin{remark}
  \label{rem_charac_sphfcn}
  Suppose $K$ and $N$ are Lie groups and $G=KG^o$,
  where $G^o$ is the connected component of the neutral element.
  Then the spherical functions $\Phi$ on $G$ are analytic 
  and they are the common eigenfunction of 
  ($G$-)left and $K$-invariant differential operators on $G$,
  such that $\Phi(0)=1$.
  Equivalently, the spherical function $\phi$ on $N$
  are analytic and they are the common eigenfunction 
  of  ($N$-)left and $K$-invariant differential operators on $N$,
  such that $\phi(0)=1$  \cite[ch.X]{helgason}.
\end{remark}
As examples of spherical functions,
we shall provide their expressions on Heisenberg groups.
These will be used during the proof of Theorem~\ref{mainthm}.\\
If $C^\natural(N)$ is a commutative algebra for the convolution product,
then $(N,K)$ is called a Gelfand pair.\\
We recall the link between bounded spherical functions and representations,
which we will use to construct our bounded spherical functions:
\begin{theorem}
  \label{thm_sphfcn_rep}
  Let  $(N,K)$ be a Gelfand pair.

  \begin{itemize}
  \item[a)]\cite[ch.X]{helgason}, \cite[ch.IV,I]{faraut}
    The vector space of $K$-invariant vectors 
    of an irreducible representation on $G=K\ltimes N$
    is of dimension at most one.

    The spherical functions of positive type (on $G$)
    are the positive definite functions $\Phi$ (on $G$)
    which are associated to an irreducible representation
    with at least one non zero $K$-invariant vector.

    For the representation associated to a positive definite function,
    the vector space of $K$-invariant vectors is
    $\Cb \Phi$. 

  \item[b)]\cite[Corollary 8.4]{pgelf}
    If $N$ is a nilpotent Lie group, 
    then the bounded spherical functions 
    are the spherical functions of positive type.

  \end{itemize}
\end{theorem}
It is known \cite[Theorem 5.12]{pgelf} that 
$(N_p,SO_p)$ and consequently $(N_p,O_p)$ are Gelfand pairs.
Thus to obtain the bounded spherical functions of $(N_p,O_p)$,
we need to describe classes of representations of
$G:=O_p\ltimes N_p$.
In this section, we shall compute those of $N_p$ by the
orbit method (see \cite{pukanszky} or \cite{kirillov}).
We compute $\hat{G}$ using Mackey's Theorem 
\cite[ch.III B Theorem 2]{lipsman}, 
provided that we describe $\hat{N}_p/G$.\\
To describe $\hat{N_p}$, 
the classes of representations of $N$,
we shall  use the orbit method (see \cite{pukanszky} or \cite{kirillov}).
For a connected simply connected nilpotent Lie group $N$,
we will denote by $T_f$ the classes of representation of $N$ associated to $f\in \Nc^*$.
First, we set the following conventions for elements of $\Nc^*$.

\paragraph{Conventions regarding elements of \mathversion{bold}{$\Nc^*$}.}
In this section and in the rest of this paper, 
we write $N=N_p$,
its Lie algebra $\Nc_p=\Nc$
and the dual $\Nc_p^*=\Nc^*$.
We denote by $\Vc^*$ and $\Zc^*$
the dual spaces of $\Vc$ and $\Zc$ respectively,
and by $X_1^*,\ldots,X_p^*$ the dual basis of $X_1,\ldots,X_p$.\\
Let $A^*\in\Zc^*$ be identified with an antisymmetric transformation
(by the scalar product on $\Zc$).
We associate to it the bilinear antisymmetric form $\omega_{A^*}$
on $\Vc$, given by: 
$\omega_{A^*}(X,Y)=<A^*X,Y>$,
$X,Y\in\Vc$.
The radical of $\omega_{A^*}$ coincides with $\ker A^*$;
and its orthogonal complement in $(\Vc,<,>)$ is $\Im A^*$,
the range of $A^*$. 
So on $\Im A^*$, $\omega_{A^*}$ induces a simplectic form
$\omega_{A^*,r}$ and the dimension of $\Im A^*$ is even 
and will be denoted by $2p_0$.\\
Suppose we have fixed $E_1$, 
a maximal totally isotropic space for  $\omega_{A^*,r}$. 
Then $E_2=A^*E_1$ is 
the orthogonal complement of $E_1$ in $(\Im A^*,<,>)$ 
and  a maximal totally isotropic space for $\omega_{A^*,r}$.
The dimension of $E_1$ and $E_2$ is $p_0$.
We denote by
$q_0:\Vc\rightarrow \ker A^*$,
$q_1:\Vc\rightarrow E_1$
and $q_2:\Vc\rightarrow E_2$
the orthogonal projections.

\paragraph{Description of \mathversion{bold}{$\hat{N_p}$}.}
Now we describe $\hat{N_p}$, 
the classes of representations of $N$
(we will only need some of these classes).\\
We need to describe first the representatives of $\Nc^*/N$.
The co-adjoint representation is given for $n=\exp (X+A)\in N$ by:
$$
X^*\in\Vc^*,\; A^*\in\Zc^*
\quad
\Coad.n(X^*+A^*) 
\,=\,
X^*+A^*
-A^*.X
\quad.
$$
We can thus choose the privileged representative $X^*+A^*$
($X^*\in \Vc$ and $A^*\in\Zc$),
of each orbit $\Nc^*/N$, 
such that $X^*\in \ker A^*$. 
Let $f=X^*+A^*$ have this form.
We define the bilinear antisymmetric form on $\Nc$ associated to $f$:
$$
\forall \; V,V'\in\Nc
\quad :\quad 
B_f(V,V')\,=\,f([V,V'])
\quad.
$$
Because of~(\ref{equality_bracket_scpdt}), we have:
$$
B_f(X+A,X'+A')
\,=\,
f([X,X'])
\,=\, 
<A^*,[X,X']>
\,=\, 
w_{A^*}(X,X')
\quad .  
$$
Some easy computations show that 
a polarization $\Lc_f$ at $f$ 
and an associated representation
$(\Hc_{X^*,A^*},U_{X^*,A^*})$ are given by:
\begin{itemize}
\item if $A^*=0$, 
  then $B_f=0$, $\Lc_f=\Nc$, and
  $U_{X^*,A^*}$ is the one dimensional representation given by:
  $\exp(X+A) \mapsto\, \exp(i<X^*,X>)$.
\item if $A^*\not=0$
  (with the previous conventions about $A^*\in\Zc^*$),
  we assume that we have chosen 
  a maximal totally isotropic space $E_1$ for  $\omega_{A^*,r}$,
  and so $\Lc_f:=E_2\oplus \ker A^*\oplus\Zc$
  is a polarization at $f$.
  Another choice for $E_1$ gives another polarization at $f$,
  but does not change the class of $U_{X^*,A^*}$.
  We compute
  $\Hc_{X^*,A^*}=L^2(E_1)$,
  and for
  $F\in \Hc_{X^*,A^*}$, 
  $n=\exp(X+A)$,
  $X'\in  E_1$:
  \begin{eqnarray}
    U_{X^*,A^*} (n).F (X')
    &=&
    \exp\left(i<A^*, \frac12[q_1(X+2X'),q_2(X)]+A>\right)\nonumber\\
    && \quad e^{ i<X^*,X>}
    F\left(q_1(X)+X'\right)\quad.\label{def_UX*A*}
  \end{eqnarray}
\end{itemize}
Kirillov's Theorem gives:
\begin{proposition}
  For $A^*\in\Zc^*$, and $X^*\in\ker A^*\subset\Vc^*$,
  we have:
 $$
U_{X^*,A^*}\in T_{X^*+A^*}
\quad.
$$
  Furthermore, 
  when $A^*$ and $X^*$ ranges over $\Zc^*$ and $\ker A^*$ respectively,
  $U_{X^*,A^*}$ ranges over a set of representatives of each class of $\hat{N}_p$.
\end{proposition}
\begin{remark}  
  \label{rem_expression_ker}
  The Lie algebra of $\ker U_{X^*,A^*}$ is:
  $$
  \left( \ker A^*\cap {(X^*)}^\perp\right) \oplus {(A^*)}^\perp
  \quad,  
  $$
  where ${(X^*)}^\perp$ is the orthogonal space of $X^*$ in $(\Vc,<,>)$,
  and  ${(A^*)}^\perp$ is the orthogonal space of $A^*$ in $(\Zc,<,>)$.
\end{remark}
\begin{remark} 
  \label{rem_expression_centre}
  The restriction of  $U_{X^*,A^*}$ on $\Zc$ is given by:
  $$
  \exp A\mapsto\exp (i<A^*,A>)
  \quad.
$$
\end{remark}

\paragraph{Consequences  of Kirillov's Theorem.}
Here we give simple consequences of Kirillov's Theorem,
which will permit us to describe $\hat{N}_p/G$ (where $G=K\ltimes N_p$).\\
In this paragraph, $N$ is a connected simply connected nilpotent Lie group,
and $G$  a group which acts continuously by automorphisms
on $N$. We denote by $\Nc$ the Lie algebra of $N$, 
and by $\Nc^*$ the dual of $\Nc$.
Then $G$ acts on $\hat{N}$:
$$
g\in G,\; \rho\in\hat{N}
\qquad
g.\rho:= n\mapsto \rho(g^{-1}.n)
\;,
$$
and by automorphisms on the vector space $\Nc^*$:
$$
g\in G,\; f\in\Nc^*\qquad
g.f:= n\mapsto f(g^{-1}.n)
\;.
$$
For $g\in G$, we compute:
$ g.T_f=T_{g.f}$.
We deduce:
\begin{corollary}
  \label{cor_thm_kirillov}
  The Kirillov map induces a one-to-one map from $(\Nc^*/N)/G$ onto $\hat{N}/G$,
  which maps the $G$-orbit of $f\in\Nc^*$ to the $G$-orbit of $T_f$.
\end{corollary}
Under the previous hypothesis,
for $\rho\in\hat{N}$,
we denote its $G$-stability group by: 
$$
G_\rho
\,=\,  \{ g\in G;\; g.\rho= \rho \}
\quad.
$$
By Kirillov's orbit method, it is easy to see that:
\begin{proposition}
  \label{prop_stability_group}
  Let $N$ be a connected simply connected nilpotent Lie group,
  and $K$ a group which acts continuously by automorphisms
  on~$N$. Let $G=K\ltimes N$.\\
  We have $(\Nc^*/N)/G\sim\Nc^*/G$.\\
  Furthermore, let $\rho\in \hat{N}$ be fixed.
  We may assume that $\rho=T_f, f\in  \Nc^*$.
  Then the $G$-stability group $G_\rho$ 
  is $K_\rho\ltimes N$, where
  $K_\rho$ is the $K$-stability group of $\rho$,
  or equivalently of the $N$-orbit $N.f$ of $f$:
  $$
  \mbox{i.e.}\qquad
  K_\rho
  \,:=\,
  \{ k\in K:\;
  k.\rho = \rho
  \} 
  \,=\,
  \{ k\in K\subset G:\;
  k.f \in N.f
  \}  \quad.
  $$
\end{proposition}

\paragraph{Bounded Spherical Function on the Heisenberg Group.}
Here, as example of spherical functions,
we provide the expressions of the bounded spherical functions
on Heisenberg groups for some compact groups. 
This will be used during the proof of Theorem~\ref{mainthm}.\\
We use the following law of 
the Heisenberg group $\Hb^{p_0}$:
\begin{eqnarray*}
  \forall \, h=(z_1,\ldots,z_{p_0},t)\, ,\,
  h'=(z'_1,\ldots,z'_{p_0},t') \in \Hb^{p_0}=\Cb^{p_0}\times \Rb  \\
  h.h'=
  (z_1+z'_1,\ldots,z_{p_0}+z'_{p_0}  ,   
  t+t'+\frac12\sum_{i=1}^{p_0}  \Im  z_i\bar{z}'_i)
  \quad.
\end{eqnarray*}
The unitary $p_0\times p_0$ matrix group $U_{p_0}$
acts by automorphisms on $\Hb^{p_0}$.
Let us describe some subgroups of $U_{p_0}$.
Let $p_0,p_1\in \Nb$,
and $m=(m_1,\ldots,m_{p_1})\in \Nb^{p_1}$ 
be fixed
such that $\sum_{j=1}^{p_1} m_j =p_0$.
We define $m'_j$ for
$j=1,\ldots p_1$ by~(\ref{formula_m_m'}).
Let $K(m;p_1;p_0)$ be the subgroup of~$U_{p_0}$ given by:
\begin{equation}
  \label{K(m,p)}
  K(m;p_1;p_0)=U_{m_1}\times \ldots \times U_{m_{p_1}}
  \quad.
\end{equation}
The expressions of spherical functions of $(\Hb^{p_0},K(m;p_1;p_0))$
can be found in the same way as in the case 
$m=(p_0), p_1=1$ i.e. $K=U_{p_0}$ \cite[ch.V,II.6]{faraut_h}
using Remark~\ref{rem_charac_sphfcn};
here, we admit \cite{moi}:
\begin{proposition}
  \label{prop_sphfcn_heis}
  $(\Hb^{p_0},K(m;p_1;p_0))$ is a Gelfand pair.\\
  Its bounded spherical functions on $\Hb^{p_0}$ are:
  \begin{enumerate}
  \item 
    $\omega=\omega_{{\lambda},l}$ 
    with
    ${\lambda}\in \Rb^* $ 
    and $l=(l_1,\ldots,l_{p_1})\in\Nb^{p_1}$ : 
    $$
    \omega(z_1,\ldots,z_{p_0},t)
    \,=\,
    e^{i\lambda t}
    \overset{p_1}{\underset{j=1}{\Pi}}
    \flb {l_j} {m_j-1} {\frac {\nn{\lambda}}2 \sum_{m'_{j-1}<i\leq m'_j}\nn{z_i}^2} 
    \quad,
    $$
  \item 
    $\omega=\omega_\mu $
    with
    $\mu=(\mu_1,\ldots,\mu_{p_1})$ 
    and ${{\mu}}_i >0$ :
    $$
    \omega(z,t)
    \,=\,
    \overset{p_1}{\underset{j=1}{\Pi}}
    \fb {m_j-1} {{{\mu}}_j \sqrt{\sum_{m'_{j-1}<i\leq m'_j}\nn{z_i}^2}}
    \quad.
    $$
  \end{enumerate}
\end{proposition}
During the proof of Theorem~\ref{mainthm}, 
we will use the following notations.
To a spherical function~$\omega$ for the Gelfand pair~$(\Hb^{p_0},K(m;p_0;p_1))$,
we associate the corresponding spherical function $\Omega^\omega$ 
on $H_{heis}=K(m;p_0;p_1)\ltimes \Hb^{p_0}$,
and the irreducible representation $(\Hc_\omega,\Pi_\omega)$ on $H_{heis}$
associated with $\Omega^\omega$.
  We compute easily:
  \begin{eqnarray}
\mbox{if}\,\omega=\omega_{\lambda,l}\qquad  
  \Pi_\omega(0,t)&=&\exp (i\lambda t) \label{Piomega_Z_lambda}\\
\mbox{if}\,\omega=\omega_\mu \qquad  
\Pi_\omega(0,t)&=&1\quad.\label{Piomega_Z_mu}
  \end{eqnarray}


\section{Expression of the Bounded Spherical Functions}
\label{sec_proof}

This section is devoted to the proof of Theorem~\ref{mainthm} for $K=O_p$.
Let $G$ be the group $K\ltimes N$,
where $N=N_p$ and $K=O_p$.
We fix the Haar measure $dkdn$ on $G$.

\paragraph{Overview of the proof.}
For $\rho\in \hat{N}$, we denote by:
\begin{itemize}
\item $G_\rho$ the $G$-stability group of $\rho$,
\item $\check{G_\rho}$ the set of
  $\nu\in\hat{G_\rho}$ such that
  $\nu_{|N}$ is a multiple of $\rho$,
\item  $\tilde{G}_\rho$ the set of 
  $\nu\in\check{G}_\rho$  
  such that 
  the dimension of the space of $K_\rho$-invariant vectors is one.
\end{itemize}
By Mackey's Theorem  \cite[ch.III B Theorem 2]{lipsman},
when $\rho$ and $\nu$ range over a representative of each class 
of  $\hat{N}$ and $\check{G_\rho}$ respectively,
the representation  induced by~$\nu$ on~$G$
gives a representative of each class of $\hat{G}$.\\
Because of the subgroup 
and intertwining number Theorems
\cite[ch.II A, Theorem 1 and Lemma 5 respectively]{lipsman},
we easily get for $\nu\in \hat{G}_\rho$ :
\begin{eqnarray}
  \nu\in\tilde{G}_\rho
  &\Longleftrightarrow&
  \nu_{|K_\rho}\,
  \mbox{contains exactly one times}\,
  1_{K_\rho}
  \nonumber
  \\
  &\Longleftrightarrow&
  \mbox{the space of}\, K
  \mbox{-invariant vectors of}
  \,\Ind_{G_\rho}^G \, \nu\,
  \mbox{is a line.}
  \label{proporty_tildeGrho}
\end{eqnarray}
The proof of Theorem~\ref{mainthm} is based 
on the two theorems and proposition
which follow.
We will explain after their statements 
how we deduce from them the expression of all bounded spherical functions.\\
First, we express the bounded spherical functions in terms of representations
$\rho\in\hat{N}/G$ and $\nu\in\tilde{G}_\rho$:
\begin{theorem}  
  \label{thm_tildeGrho}
  Let $\rho\in \hat{N}$ 
  and $(\Hc^\nu,\nu)\in \tilde{G_\rho}$.
  Then because of (\ref{proporty_tildeGrho}),
  $\Ind_{G_\rho}^G \, \nu \,\in\, \hat{G}$
  has also a (non-zero) $K$-invariant line and
  the  associated bounded spherical function
  is the function $\phi^\nu$ given by :
  \begin{equation}
    \label{fcnsphnu}
    \phi^\nu(n)
    \,=\,
    \int_K
    {\big<\nu(I,k.n).\vec{u}_\nu\,,\,
      \vec{u}_\nu \big>}_{\Hc^\nu} dk,
    \quad 
    n\in N
    \quad ,
  \end{equation}
  where  $\vec{u}_\nu\in\Hc^\nu$ 
  is any unit $K_\rho$-invariant  vector.\\
  Furthermore, we obtain all the bounded spherical functions
  as $\phi^\nu$ 
  when $\rho$ and $\nu$
  range over a set of representatives of~$\hat{N}/G$,
  and~$\tilde{G}_\rho$ respectively.
\end{theorem}
Next, to obtain all representations $\rho\in\hat{N}/G$,
we describe $\Nc^*/G$ (see Corollary~\ref{cor_thm_kirillov}):
\begin{proposition}
  \label{prop_Nc^*/G}
Let $O(r,\Lambda)=G.(rX^*_p+D_2(\Lambda))\subset \Nc^*$.\\
Then the mapping
$$
\begin{array}{lcr}
\cM&\rightarrow&\Nc^*/G\\
(r,\Lambda)&\mapsto&O(r,\Lambda)
\end{array}
$$
is a bijection.
\end{proposition}
Now, we describe $\tilde{G}_\rho$,
where $\rho$ is a representation associated (by Kirillov) to a linear form on $\Nc$,
which is a privileged representative of a $G$-co-adjoint orbit
(just given in Proposition~\ref{prop_Nc^*/G}):
\begin{theorem}
  \label{thm_fcnsphnu}
  Let $\rho\in T_f$ where $f=rX^*_p+D_2(\Lambda)$ and $(r,\Lambda)\in \cM$.
  \begin{itemize}
  \item[a)] If $\Lambda=0$,
    $\tilde{G}_\rho\,=\,\{\nu^{r,0}\}$.
    The spherical function $\phi^\nu$
    which is  associated (by~(\ref{fcnsphnu})) to $\nu=\nu^{r,0}$ is $\phi^{r,0}$
    (given by~(\ref{def_sph_fcn_bes})).
  \item[b)] If $\Lambda\not=0$,
    $\tilde{G}_\rho    \subset    \{\nu^{r,\Lambda,l}\; ,\; l\in \Nb^{p_1}\}$.
    Each representation $\nu=\nu^{r,\Lambda,l}\in\hat{G_\rho}$ 
    has a $K_\rho$-invariant line,
    and the spherical function $\phi^\nu$
    associated (by~(\ref{fcnsphnu})) is $\phi^{r,\Lambda,l}$ (given by~(\ref{def_sph_fcn_lag})).
  \end{itemize}
\end{theorem}
The representations $\nu^{r,0}$ 
and $\nu^{r,\Lambda,l}$
will be described during the proof 
(see (\ref{formula_nu_x}) and (\ref{formula_nu_lambda})).\\
For the moment, we will admit these two theorems and the proposition,
and keep their notations.
From Corollary~\ref{cor_thm_kirillov}
and Proposition~\ref{prop_Nc^*/G},
we deduce that:
$$
\hat{N}/G
\,=\,
\{T_{rX_p^*+D_2(\Lambda)}\;,\;
(r,\Lambda)\in  \cM
\}
\quad.
$$
Under Theorems~\ref{thm_tildeGrho} and~\ref{thm_fcnsphnu},
the spherical bounded functions are 
the functions
$\phi^{r,0}$, when $r\in \Rb^+$,
and $\phi^{r,\Lambda,l}$ 
when   $(r,\Lambda)\in \cM$
and $l\in \Nb^{p_1}$.\\
If we prove
Theorems~\ref{thm_tildeGrho} and \ref{thm_fcnsphnu},
and Proposition~\ref{prop_Nc^*/G},
Theorem~\ref{mainthm} will follow.
The rest of this section will be devoted to this.
We start with the proofs of Theorem~\ref{thm_tildeGrho}
and Proposition~\ref{prop_Nc^*/G}.
Then for a representation $\rho\in T_{rX_p^*+D_2(\Lambda)}$,
we describe its $G$-stability group
and the quotient group~$\overline{N}=N/ \ker\rho$.
We finish with the proof of Theorem~\ref{thm_fcnsphnu}.

\paragraph{Set \mathversion{bold}{$\tilde{G}_\rho$}.}
The aim of this paragraph is to prove Theorem~\ref{thm_tildeGrho}.\\
Let $\rho\in\hat{N}$ be fixed.
Under Proposition~\ref{prop_stability_group}, 
the $G$-stability group of $\rho$ 
is $G_\rho=K_\rho\ltimes N$,
where $K_\rho$ is the $K$-stability group 
(which is a compact subgroup of $K$).
We fix the normalized Haar measure $dk_\rho$ on $K_\rho$,
and the Haar measure $dk_\rho dn$ on $G_\rho$.\\
We fix $(\Hc^\nu,\nu) \in \tilde{G}_\rho$
and  a unit $K_\rho$-invariant vector 
$\vec{u}=\vec{u}_\nu\in \Hc^\nu$.
We denote by $(\Hc^\Pi,\Pi)$ the induced representation $\Ind_{G_\rho}^G \, \nu$ of $\nu$:
$$
\forall\,g,g'\in G\quad
f\in \Hc^\Pi\quad:\quad
\Pi(g).f(g')\,=\, f(g'g)
\quad,
$$
and by $f$  the function on $G$ given by :
$$
f(k,n)\,=\,\nu(I,n).\vec{u}\,,\quad (k,n)\in G
\quad;
$$
the vector $f\in \Hc^\Pi$ is $K$-invariant and of norm one.
We can then associate to $\Pi$ and $f$ the bounded spherical function $\phi^\nu$:
$$
\phi^\nu(g)
\,=\,
{\big<\Pi(g).f,f\big>}_{\Hc^\Pi}
\,,\quad g\in G\quad.
$$
We can easily obtain for $g=(k,n),g'=(k',n')\in G$:
\begin{eqnarray*}
  \Pi(g).f(g')
  &=&
  \nu(I,n')\nu(I,k'.n).\vec{u}
  \quad,
  \\
  {\big< \Pi(g).f(g'),f(g') \big>}_{\Hc^\nu}
  &=&
  {\big<\nu(I,k'.n).\vec{u}\,,\,
    \vec{u} \big>}_{\Hc^\nu}
  \quad .
\end{eqnarray*}
We thus obtain the formula~(\ref{fcnsphnu}).\\
We can now complete our proof of Theorem~\ref{thm_tildeGrho}.
Under Mackey's Theorem and property~(\ref{proporty_tildeGrho}), 
when $\rho$ and $\nu$ range over a set of representatives of $\hat{N}/G$ 
and $\tilde{G}_\rho$ respectively, 
we get all the irreducible representations
$\Pi=\Ind_{G_\rho}^G \, \nu$ having a $K$-invariant line. 
Under Theorem~\ref{thm_sphfcn_rep},
the positive definite  functions $\phi^\nu$ associated to $\Pi$
give all the bounded spherical functions.

\paragraph{Description of \mathversion{bold}{$\Nc^*/G$}.}
Here we prove Proposition~\ref{prop_Nc^*/G}.
We easily compute
for $g=(k,n)\in G$ with $n=\exp (X+A)\in N$
and $X^*\in\Vc^*,A^*\in\Zc^*$:
\begin{equation}
  \label{equality_coad}
  \Coad.g(X^*+A^* )
  \,=\,
  k.X^*+k.A^*-(k.A^*).X
  \quad . 
\end{equation}
Let $O\in \Nc^*/G$ be a fixed orbit.
We associate to it $\Lambda\in\overline{\Lc}$
such that all the antisymmetric matrices $A^*_f$,
where $f=X^*_f+A^*_f\in O$, are $K$-conjugate,
and $K$-conjugated to  $D_2(\Lambda)$.\\
Let $f=X^*_f+A^*_f\in O$ be fixed. 
We make the following choices:
\begin{enumerate}
\item let $k_0\in K$ be such that  $k_0.A_f^*=D_2(\Lambda)$; 
\item let $X_0\in \Vc$ be 
  such that $(k_0.A_f^*).X_0\in \Vc^*$ is the orthogonal projection $X^*_0$ of $k_0.X_f^*\in \Vc^*$ 
  on the kernel $\ker k_0.A^*_f=\ker D_2(\Lambda)$;
  in particular, $X^*_0=0$ if $\Im  D_2(\Lambda)=\Vc$;
\item let $k'_0\in K$ be such that $k'_0.X\in \Im D_2(\Lambda)$
  for all $X\in \Im D_2(\Lambda)$
  and $k'_0X_0^*=rX_p^*$, $r\in\Rb^+$.
\end{enumerate}
We get $(k_0'k_0,\exp X_0).f=rX^*_p+D_2(\Lambda)$.\\
We remark that $\Im D_2(\Lambda)=\Vc$ is equivalent to $p=2p_0$
and $\lambda_i\not=0$, $i=1,\ldots,p'$.\\
Proposition~\ref{prop_Nc^*/G} is thus proved.

\paragraph{Stability Group \mathversion{bold}{$K_\rho$}.}
The aim of this paragraph is to describe the stability group $K_\rho$
of $\rho\in T_{rX^*_p+D_2(\Lambda)}$.\\
Before this, 
let us recall that the orthogonal $2n\times 2n$ matrices which commutes with
$D_2(1,\ldots,1)$ (see (\ref{def_D2Lambda}) for this notation) 
have determinant one and form the group  $Sp_n\cap O_n$.
This group is isomorphic to $U_n$;
the isomorphism is denoted $\psi_1^{(n)}$,
and satisfies:
$$  \forall\,k,X
\quad:\quad 
\psi_c^{(n)} (k.X) 
\,=\,
\psi_1^{(n)}(k).\psi_c^{(n)}(X)
\quad,  
$$
where $\psi_c^{(n)}$ is the complexification:
$$
\psi_c^{(n)}(  x_1,y_1;\ldots;x_n,y_n)=
(x_1+i y_1,  \ldots,  x_n+iy_n)
\quad.
$$
Now, we can describe $K_\rho$:
\begin{proposition}
  \label{prop_Krho}
  Let $(r,\Lambda)\in \cM$. 
  Let $p_0$ be the number of $\lambda_i\not=0$, 
  where $\Lambda=(\lambda_1,\ldots,\lambda_{p'})$,
  and $p_1$ the number of distinct $\lambda_i\not=0$.
  We set
  $\tilde{\Lambda}=(\lambda_1,\ldots,\lambda_{p_0})\in \Rb^{p_0}$.

  Let $\rho\in T_f$ where $f=rX^*_p+D_2(\Lambda)$.
  \begin{itemize}
  \item If $\Lambda=0$, then $K_\rho$ is the subgroup of~$K$
    such that  $k.rX^*_p=rX^*_p$ for all $k\in K_\rho$. 
  \item If $\Lambda\not=0$, 
    then  $K_\rho$ is the direct product $K_1\times K_2$, 
    where:
    \begin{eqnarray*}
      K_1
      &=&
      \left\{
        k_1 \,=\,
        \left[ \begin{array}{cc}
            \tilde{k}_1&0\\
            0& \Id
          \end{array}\right]
        \quad \left/\quad
          \begin{array}{l}
            \tilde{k}_1\in SO(2p_0)\\
            D_2(\tilde{\Lambda})
            \tilde{k}_1
            =
            \tilde{k}_1
            D_2(\tilde{\Lambda})
          \end{array}
        \right.\right\}
      \quad,\\
      K_2
      &=&
      \left\{
        k_2 \,=\,
        \left[ \begin{array}{cc}
            \Id&0\\
            0& \tilde{k}_2
          \end{array}\right]\in K
        \quad\left/\quad
          \tilde{k}_2.rX^*_p=rX^*_p
        \right.\right\}
      \quad.
    \end{eqnarray*}
    Furthermore,
    $K_1$ is isomorphic to the group $K(m;p_0;p_1)$ given by (\ref{K(m,p)}). 
  \end{itemize}
\end{proposition}

\begin{proof}
  We keep the notations of this proposition,
  and we set $A^*=D_2 (\Lambda)$
  and $X^*=rX^*_p$.
  With Propositions~\ref{prop_stability_group} 
  and the expression~(\ref{equality_coad}) of the co-adjoint representation,
  it is easy to prove:
  \begin{equation}
    \label{formula_Krho}
    K_\rho
    \,=\, \{ k\in K\; :\; 
    kA^*=A^*k
    \quad\mbox{and}\quad
    kX^*=X^* \}
    \quad.
  \end{equation}
  If $\Lambda=0$, because of~(\ref{formula_Krho}),
  $K_\rho$ is the stability group in $K$ 
  of $X^*\in \Vc^*\sim \Rb^p$.
  So the first part of Proposition~\ref{prop_Krho} 
  is proved.

  Let us show the second part.
  $\Lambda\not =0$ so we have
  $$
  A^*
  \,=\,
  \left[\begin{array}{c|c}
      D_2(\tilde{\Lambda})&0\\
      \hline  0&0
    \end{array}\right]
  \quad\mbox{with}\quad
  D_2(\tilde{\Lambda})
  \,=\,
  \left[\begin{array}{ccc}
      \mu_1 J_{m_1}&0&0\\
      &\ddots&\\
      0&0&\mu_{p_1} J_{m_{p_1}}
    \end{array}\right]
  \quad.
  $$
  where $\mu_1,\ldots,\mu_{p_1}$ are defined by (\ref{formula_mu_lambda}),
  and $m_j$ is the number of $\lambda_i=\mu_j$.
  We define $m'_j$ for
  $j=1,\ldots p_1$ by~(\ref{formula_m_m'}).\\
  Let $k\in K_\rho$. 
  Because of~(\ref{formula_Krho}), the matrices $k$ and $A^*$ commute and we have:
  $$
  k \,=\,
  \left[ \begin{array}{cc}
      \tilde{k}_1&0\\
      0& \tilde{k}_2
    \end{array}\right]
  \quad\mbox{with}\; \tilde{k}_1\in O(2p_0)
  \;\mbox{and}\;
  \tilde{k}_2\in O(p-2p_0)
  \quad;
  $$
  furthermore,
  by~(\ref{formula_Krho}),
  $\tilde{k}_2.X^*=X^*$,
  and the matrices $\tilde{k}_1$ and $D_2(\tilde{\Lambda}^*)$  commute.
  So $\tilde{k}_1$ is a diagonal block matrix, with blocks
  ${[\tilde{k}_1]}_j \in O(m_j)$ for $i=1,\ldots ,p_1$.
  Each block ${[\tilde{k}_1]}_j \in O(m_j)$ 
  commutes with  $J_{m_j}$.
  So on one hand we have 
  $\det {[\tilde{k}_1]}_j =1$,  
  $\det\tilde{k}_1=1$,
  and on the other hand, 
  ${[\tilde{k}_1]}_j \in O(m_j)$ 
  corresponds to a unitary matrix $\psi_1^{(m_j)}({[\tilde{k}_1]}_j)$.
  Now, we set for $k_1\in K_1$:
  $$
  \Psi_1(k_1)
  \,=\,
  \left( \psi_1^{(m_1)}({[\tilde{k}_1]}_1),\ldots,
    \psi_1^{(m_{p_1}}({[\tilde{k}_1]}_{p_1})\right)
  \quad .  
  $$
  $\Psi_1 : K_1\rightarrow K(m;p_0;p_1)$ 
  is a group isomorphism.
\end{proof}

As $\psi_1^{(n)}$ is an isomorphism which respects complexification, 
we have:
\begin{corollary}
  \label{corollary_Krho}
  The isomorphism $\Psi_1$ given during the previous proof
  respects complexification:
  $$
  \psi_c^{(p_0)}\{ \tilde{k}_1.(x_1,y_1,\ldots,x_{p_0},y_{p_0})\}
  \,=\,
  \Psi_1(k_1)
  .\psi_c^{(p_0)}(x_1,y_1,\ldots,x_{p_0},y_{p_0})
  \quad.
  $$
\end{corollary}

\paragraph{Quotient Group \mathversion{bold}{$\overline{N}=N/\ker\rho$}.}
In this paragraph, 
we describe the quotient groups
$N/\ker \rho$ and $G/\ker \rho$, for some $\rho\in\hat{N}$.
This will permit in the next paragraph
to reduce the construction of the bounded spherical functions on $N_p$
to known questions on Euclidean and Heisenberg groups.
For a representation
$\rho\in\hat{N}$,
we will denote by:
\begin{itemize}
\item $\ker\rho$ the kernel of $\rho$,
\item $\overline{N}=N/\ker\rho$ its quotient group
  and $\overline{\Nc}$ its Lie algebra,
\item  $(\Hc,\overline{\rho})$ the induced representation on $\overline{N}$,
\item $\overline{n}\in\overline{N}$ and $\overline{Y}\in\overline{\Nc}$ 
  the image of $n\in N$ and $Y\in\Nc$ respectively
  by the canonical projections
  $N \rightarrow \overline{N}$ and
  $\Nc \rightarrow \overline{\Nc}$.
\end{itemize}
Now, with the help of the canonical basis,
we choose the privileged representative $\rho$
of $T_{rX^*_p+D_2(\Lambda)}$
as $\rho=U_{rX^*_p,D_2(\Lambda)}$ given by~(\ref{def_UX*A*}), 
with $E_1=\Rb X_1\oplus\ldots\oplus\Rb X_{2p_0-1}$
as maximal totally isotropic space for  $\omega_{D_2(\Lambda),r}$. 
Because of Remark~\ref{rem_expression_ker},
the quotient Lie algebra $\overline{\Nc}$ has the natural basis:
$$
\overline{X_1},\ldots,\overline{X_{2p_0}},  
\overline{B}=\nn{\Lambda}^{-1} \overline{D_2(\Lambda)} 
\quad
\mbox{with}\;X_p\;\mbox{if}\;r\not=0\quad;
$$
here, we have denoted 
$\nn{\Lambda}={(\sum_{j=1}^{p'} {\lambda_j}^2)}^\frac12=\nn{D_2(\Lambda)}$
(for the Euclidean norm on $\Zc$). 
We compute that each Lie bracket of two vectors of this basis equals zero,
except:
$$
[\overline{X_{2i-1}},\overline{X_{2i}}] 
\,=\,
\frac{\lambda_i }{\nn{\Lambda}}\overline{B},
\quad i=1,\ldots,p_0
\quad.
$$
Let $\overline{\Nc_1}$ be the Lie sub-algebra of $\overline{\Nc}$,
with basis
$\overline{X_1},\ldots,\overline{X_{2p_0}},  \overline{B}$,
and let $\overline{N_1}$ be 
its corresponding connected simply connected nilpotent Lie group.
We define the mapping $\Psi_2 : \Hb^{p_0}\rightarrow\overline{N_1}$
for $h=(x_1+iy_1,\ldots, x_{p_0}+iy_{p_0},t)\in\Hb^{p_0}$ by:
$$
\Psi_2(h)
\,=\,
\exp \left(\sum_{j=1}^{p_0}  
  \sqrt{\frac{\nn{\Lambda}}{\lambda_j}}
  \left(x_j\overline{X_{2j-1}}+y_i\overline{X_{2j}}\right)
  \;
  +t\overline{B}\right)
\quad.
$$
Because  of the values of the Lie brackets in $\overline{N_1}$, 
it is easy to see:
\begin{lemma}
  \label{lem_isom_barN1_heis}
  $\Psi_2$ is a group isomorphism between  $N_1$ and   $\Hb^{p_0}$.   
\end{lemma}
With our choice for representations and notations, 
we describe the induced representation and action: 
\begin{proposition}
  \label{prop_barN_barrho}
  With the notations set just above,
  \begin{itemize}
  \item[a)] If $\Lambda\not=0$,
    the two groups $\overline{N}$ and 
    $\overline{N_1}\times\overline{N_2}$ are isomorphic
    and the representations
    $\overline{\rho}$ and $ \overline{\rho_1}\otimes \overline{\rho_2}$ 
    are equivalent,
    where
    \begin{enumerate}
    \item $\overline{\rho}_1$ is a representation on $\overline{N_1}$ 
      (whose expression may be computed);
    \item $\overline{N_2}$ and $\overline{\rho}_2$ are described by:
      \begin{itemize}
      \item either $r=0$, 
        then $\overline{N_2}$ and $\overline{\rho}_2$ are trivial;
      \item or $r\not =0$, 
        then $\overline{N_2}\sim\Rb \overline{X}_p$,
        and $\overline{\rho}_2\sim
        \exp (x \overline{X}_p)  \rightarrow \exp (ix)$.
      \end{itemize}
    \end{enumerate}
  \item[b)] If $\Lambda=0$,
    then $\overline{N}$ and $\overline{\rho}$ 
    are the same as 
    $\overline{N_2}$ and $\overline{\rho}_2$ above.
  \end{itemize}
\end{proposition}
\begin{remark}\label{rem_expression_centre2}
  Because of Remark~\ref{rem_expression_centre},
  the restriction of $\overline{\rho}_1$ on the center $\exp\Rb \overline{B}$ of $\overline{\Nc_1}$ 
  is given by: $\exp (a \overline{B})\mapsto \exp (ia\nn{\Lambda})$.
\end{remark}
As $K_\rho$ is  the $K$-stability group of $\rho\in\hat{N}$,
it acts by automorphisms on $\overline{N}$. 
Simple computations show that:
\begin{proposition}
  \label{prop_Krho_barN}
  We keep the notations of Propositions~\ref{prop_Krho}
  and~\ref{prop_barN_barrho}.
  \begin{itemize}
  \item[a)]  If $\Lambda=0$,
    $K_\rho$ acts trivially on $\overline{N}$.
  \item[b)] If $\Lambda\not=0$,
    \begin{itemize}
    \item $K_1$ acts by automorphisms on $\overline{N_1}$ 
      and trivially on the center $\exp\Rb \overline{B}$ of  $\overline{N_1}$, 
    \item $K_2$ acts trivially on~$\overline{N_2}$.
    \end{itemize}
    So the groups $K_\rho \ltimes\overline{N}$
    and $(K_1\ltimes \overline{N_1} )
    \,\times\, K_2\,\times\, \overline{N_2}$ are isomorphic.
  \end{itemize}
\end{proposition}
Recall $H_{heis}= K(m;p_0;p_1) \ltimes \Hb^{p_0}$;
let us also define the group 
$H=K_1\ltimes \overline{N_1}$
and the map:
$$
\Psi_0\;:\;\left\{
  \begin{array}{rcl}
    H_{heis}
    &\longrightarrow& 
    H\\
    (k_1\, , \,
    h)
    &\longmapsto&
    (\Psi_1(k_1)\, ,\,
    \Psi_2 (h))
  \end{array}\right.
\quad .
$$
Corollary~\ref{corollary_Krho}, Proposition~\ref{prop_Krho_barN},
and Lemma~\ref{lem_isom_barN1_heis}
imply:
\begin{proposition}
  \label{prop_isom_barN1_heis}
  $\Psi_0$ is a group isomorphism between 
  $H_{heis}$ and   $H$. 
\end{proposition}
\paragraph{Expression of \mathversion{bold}{$\phi^\nu$}.}
Here we prove Theorem \ref{thm_fcnsphnu}.
Let $\rho\in\hat{N}$ 
and $\nu\in\tilde{G}_\rho$ be fixed.
We have $\nu_{|N}=c.\rho$, $1\leq c \leq \infty$,
and we denote by
$\bar{\nu}$ the induced representation 
on $K_\rho\ltimes\overline{N}$.
\paragraph{a) Case of the orbit \mathversion{bold}{$O(r,0)$}.}
We assumed that $\rho=\rho_{r,0}$.
By Proposition~\ref{prop_Krho_barN},
we have $K_\rho\ltimes\overline{N}=K_\rho\times\overline{N}$.
So $\overline{\nu}$ is the tensor product of
an irreducible  representation over $\overline{N}$, 
which coincides with $c.\overline{\rho}$
(and so $c=1$),
with an irreducible representation over $K_\rho$
with a $K_\rho$-invariant vector, 
which is thus the trivial representation.
We obtain that $\overline{\nu}$ coincides
with
$(k,n)\in K_\rho\times \overline{N}\mapsto \overline{\rho}(n)$.
Now because of Proposition~\ref{prop_barN_barrho}.b),
$\nu=\nu^{r,0}$ is given by:
\begin{equation}
  \label{formula_nu_x}
  \left(k,\exp (X + A) \right)
  \,\longmapsto\,
  e^{i <rX^*_p,X>}
  \quad.
\end{equation}
So $\tilde{G}_\rho$ is the set of the classes of $\nu^{r,0}$,
where  $r$ ranges over $\Rb^+$,
and we compute that the function
$\phi^\nu$ for $\nu=\nu^{r,0}$ is given by~(\ref{fcnsphnu}):
$$
\phi^\nu(n)
\,=\,
\int_{k\in K}
e^{i <rX^*_p,k.X>}
dk,
\quad n=\exp (X+A)\in N_p
\quad.
$$  
By~(\ref{formula_bessel}),
we have
$\phi^\nu=\phi^{r,0}$,
and Theorem~\ref{thm_fcnsphnu}.a) is proved.
\paragraph{b) Case of the orbit \mathversion{bold}{$O(r,\Lambda)$}.}
We assume $\Lambda\not=0$ and $\rho=\rho_{r,\Lambda}$.
For each bounded spherical function $\omega$ 
of the Gelfand pair~$(\Hb^{p_0},K(m;p_0;p_1))$,
we define  the representation
$(\Hc^\omega,\Pi^\omega)$ 
of~$H$ such that :
$$
\Hc^\omega=
\{F\circ \Psi_0^{-1}\; , \quad F\in \Hc_\omega\}
\quad\mbox{and}\quad
\Pi^\omega=\Pi_\omega\circ\Psi_0^{-1}
\quad.
$$
By Proposition~\ref{prop_Krho_barN}.b),
$\overline{\nu}$ is the tensor product of three 
irreducible representations,
of $\overline{N_2}$,
$H$ and $K_2$,
such that 
the vector space of $K_1$ (respectively $K_2$)-invariant vectors 
of the representations of $H$ (respectively $K_2$) is a line.
So the representation of $K_2$ is trivial,
and $\overline{\nu}$ induces a unitary irreducible representation $\overline{\overline{\nu}}$ 
of~$H\times  \overline{N_2}$
which coincides with~$c.\overline{\rho}$ on $\overline{N}$,
and  such that 
the vector space of $K_1$-invariant vectors is a line.\\
We have  $\overline{\overline{\nu}}=\gamma_1\otimes\gamma_2$ where
\begin{itemize}
\item[(a)] 
  $\gamma_1$ is an irreducible representation of $H$;
  the space of its  $K_1$-invariant vectors is a line;
  it  coincides with
  $c.\overline{\rho}$ over $\overline{N_1}$;
\item[(b)]
  $\gamma_2$ is an irreducible representation over $\overline{N_2}$
  and coincides with
  $c.\overline{\rho}$ on $\overline{N_2}$.
\end{itemize}
Because of irreducibility, (a) implies $c=1$ and so (b) implies $\gamma_2\sim \rho_2$.
Furthermore,
by Proposition~\ref{prop_Krho_barN}
and Theorem~\ref{thm_sphfcn_rep},
the irreducible representations on~$H$ 
such that the vector space of $K_1$-invariant vectors is a line,
are all the representations~$(\Hc_\omega,\Pi_\omega)$,
where $\omega$ ranges over the set of bounded spherical functions of~$H_{heis}$;
the $K_1$-invariant line of  $\Hc^\omega$ is
$\Cb \Omega^\omega\circ\Psi_0^{-1}$.
Thus, the representations~$\gamma_1$ 
satisfying~(a) are the representations such that
$\gamma_1\sim\Pi^\omega$ and
$\Pi^\omega_{|\overline{N_1}}\sim\overline{\rho}_1$.
Because of the expressions of $\gamma_1$ and $\Pi^\omega$
on the center of $\overline{N_1}$ and $H$
(see Remark~\ref{rem_expression_centre2}
and equalities
(\ref{Piomega_Z_lambda}),
(\ref{Piomega_Z_mu}))
the case $\omega=\omega_\mu$ is impossible 
if  $\Pi^\omega_{|\overline{N_1}}\sim\overline{\rho}_1$.\\
We have shown that
$\overline{\overline{\nu}}=\gamma_1\otimes\gamma_2$,
where $\gamma_2\sim\rho_2$ is given in Proposition~\ref{prop_barN_barrho},
and $\gamma_1$ is among the representations equivalent to
$(\Hc_\omega,\Pi_\omega)$ with $\omega=\omega_{\lambda,l}$, $l\in\Nb^{p_1}$.\\
Then $\nu$ is  among the representations equivalent to $(\Hc^\omega,\nu^{r,\Lambda,l})$
with  $\omega=\omega_{\nn{\Lambda},l}$,
defined for $    n=\exp (X+A)\in N$, 
and $k=k_1k_2\in K_\rho$ where $k_1\in K_1$, and $k_2\in K_2$ by:
\begin{equation}
  \label{formula_nu_lambda}
  \nu^{r,\Lambda,l} (k,n) 
  \,=\,
  e^{i r<X^*_p,X>}
  \, \Pi^\omega (k_1, \bar{q_1}(n))
  \quad,  
\end{equation}
where  $\bar{q_1}:N\rightarrow \overline{N_1}$ is
the canonical projection.\\
We denote by $\tilde{G}_\rho '$ the set of classes
of the representations~$\nu^{r,\Lambda,l}$, 
$l\in\Nb^{p_1}$.
A representation~$\nu^{r,\Lambda,l}$, 
has a unitary $K_\rho$-invariant vector
$\vec{u}=\Omega^\omega\circ \Psi_0^{-1}$.\\
We still have to show that under formula~(\ref{fcnsphnu}),
the function $\phi^\nu$ for $\nu=\nu^{r,\Lambda,l}$
satisfies $\phi^\nu=\phi^{r,\Lambda,l}$.
For $n=\exp(X+A)$, it is given  by:
\begin{equation}
  \phi^\nu(n)
  \,=\,
  \int_{K}
  e^{i r<X^*_p,k.X>}
  \omega\circ
  \Psi_2^{-1}\circ \bar{q_1} (k.n)   
  dk \quad,
  \label{egalite_phi_omega_Psi_q}
\end{equation}
where $\omega=\omega_{\nn{\Lambda},l}$. 
We compute 
\begin{eqnarray*}
  \nu (I,n) \vec{u}
  &=&
  e^{i r<X^*_p,X>}
  \, \Pi^\omega (I,
  \bar{q_1}(n))
  \Omega^\omega \circ \Psi_0
  \quad,\\
  {\big<\nu(I,n).\vec{u}\,,\,
    \vec{u} \big>}_{\Hc^\nu}
  &=&
  e^{i r<X^*_p,X>}
  \,
  {\big< \Pi_\omega (I,
    \Psi_2^{-1}\circ \bar{q_1}(n) )
    \Omega^\omega ,\, 
    \Omega^\omega \big>}_{H_\omega}
  \quad.
\end{eqnarray*}
As $\Omega^\omega$ is the positive definite function
associated to  $\Pi_\omega$, we have:
\begin{eqnarray*}
  {\big< \Pi_\omega (I,
    \Psi_2^{-1}\circ \bar{q_1} (n)   )
    \Omega^\omega ,\, 
    \Omega^\omega \big>}_{H_\omega}
  &=&
  \Omega^\omega\left( I,
    \Psi_2^{-1}\circ \bar{q_1} (n)  \right)\\
  &=&
  \omega\circ
  \Psi_2^{-1}\circ \bar{q_1} (n)   
  \quad.
\end{eqnarray*}
We know the expression $\omega=  \omega_{\nn{\Lambda},l}$ (Proposition~\ref{prop_sphfcn_heis}),
and we compute those of $\Psi_2^{-1}$ and $\bar{q_1}$.
Assuming~(\ref{formula_mu_lambda}),
we obtain:
$$
\Theta^{r,\Lambda,l}(\exp (X+A))
\,=\,
e^{i r<X^*_p,X>}
\omega\circ
\Psi_2^{-1}\circ \bar{q_1}(\exp(X+A))
\quad,
$$
where
$\Theta^{r,\Lambda,l}$ is given by~(\ref{def_theta_fcn}).
Because of~(\ref{egalite_phi_omega_Psi_q}),
$\phi^\nu=    \phi^{r,\Lambda,l}$,
and Theorem \ref{thm_fcnsphnu}.b) is consequently proved.

The proof of Theorem \ref{thm_fcnsphnu} is now over.
Theorem~\ref{mainthm} is thus proved.


\section{Representation over \mathversion{bold}{$N_p$}}
\label{sec_rq_repN}

We give here the  bounded spherical functions
in terms of representations over $N_p$ 
\cite[Theorem~G]{pgelf},
which is an equivalent way for constructing them.\\
We deduce then the eigenvalues of the sub-Laplacian for the bounded spherical functions
and the expression of the radial Plancherel measure.

\paragraph{Another expression of the bounded spherical functions \mathversion{bold}{$\phi^{r,\Lambda,l}$}.}
Let us recall a few facts contained in \cite{pgelf}.
Let $(N,K)$ be a Gelfand pair.
Let $(\Hc,\Pi)\in\hat{N}$, and $K_\Pi$ its $K$-stability group.
There exists a projective representation $W_\Pi$ of $K_\Pi$ on $\Hc$,
and an orthogonal decomposition of $\Hc=\sum V_l$ into irreducible subspaces
$W_\Pi$-invariant.
For $\zeta\in \Hc$, 
let us define the function $\phi_{\Pi,\zeta}$ by:
$$
\phi_{\Pi,\zeta}(n)
\,=\,
\int_K
<\Pi(k.n)\zeta,\zeta>dk\,,
\quad n\in N\quad.
$$
The spherical functions are the $\phi_{\Pi,\zeta}$ 
for $\zeta\in V_l$, $\nn{\zeta}=1$ and $\Pi\in \hat{N}$;
the spherical function $\phi_{\Pi,\zeta}$ is independent of 
$\zeta\in V_l$, $\nn{\zeta}=1$ and of the choice of the
representative of $\Pi$.
\begin{remark}
  \label{rem_sphfcn_repN_radfcn}
  Furthermore, for a $K$-invariant function $f$ on $N$ 
  and a $K$ left-invariant differential operator $D$ on $N$, 
  $\Pi(f)$ and $d\Pi(D)$ are  $W_\Pi$-invariant;
  they equal the identity on each $V_l$ up to a constant;
  and for $V_l$, and the corresponding spherical function $\phi$,
  the constants are respectively
  $<f,\phi>$
  and the eigenvalue of $D$ for $\phi$.
\end{remark}
Now, we give the link between the given spherical functions
and the above description for each class of $\hat{N_p}$,
which corresponds to a class $O(X^*,A^*)$.\\
Before this,
we recall definitions and properties of Hermite functions.
We denote by $H_k$ the Hermite polynomial: 
$$
H_k(s)={(-1)}^k e^{s^2}{(d/ds)}^k e^{-s^2}
\quad,
$$
by $h_k,k\in\Nb$ the Hermite functions on $\Rb$ \cite[\S 5.5]{szego}:
$$
h_k(x)
\,=\,
{(2^k k!\sqrt{\pi})}^{-\frac k2} e^{-\frac{x^2}2} H_k(x)
\quad,
$$
and by $h_\alpha,\alpha\in\Nb^n$ 
the  Hermite functions on $\Rb^n$:
$h_\alpha
\,=\,
\Pi_{i=1}^n h_{\alpha_i}$.
We recall that
the Hermite functions 
$h_k,k\in\Nb$ on $\Rb$
form an orthonormal basis $L^2(\Rb)$.
\paragraph{Case $A^*=0$.}
$\rho=\rho^{r^*,0}$ is the one-dimensional representation 
given by 
$$\exp(X+A) \mapsto\, \exp(i<r X^*_p,X>)
\quad.
$$
\paragraph{Case $A^*\not=0$.}
Let $(r,\Lambda)\in\cM$
with $\Lambda\not =0$ et $l\in \Nb^{p_1}$.
Let $E_l$ be the set of 
$\alpha=(\alpha^1,\ldots,\alpha^{p_1})$ 
where $\alpha^j={(\alpha^j_i)}_{m'_{j-1}<i\leq m'_j} \in \Nb^{m_j}$
such that
$$
\nn{\alpha^j}=\sum_{m'_{j-1}<i\leq m'_j} \alpha_i^j =l_j
\quad \mbox{for}\; j=1,\ldots,p_1\quad.
$$
Let us define the representation
$(\Hc,\Pi)=(L^2(\Rb^{p_0}),\Pi_{r,\Lambda})\in\hat{N}_p$ 
for  $f\in \Hc$,
$(y_1,\ldots,y_{p_0})\in\Rb^{p_0}$,
$n=\exp(X+A)\in N_p$
by:
\begin{eqnarray*}
  \Pi(n).f(y)
  \,=\,
  e^{i<D_2(\Lambda),A> + i<r X^*_p,X>
    +i \sum_{j=1}^{p'}
    \frac {\lambda_j}2 x_{2j}x_{2j-1} +\sqrt{\lambda_j} x_{2j}y_j}\\
  f(y_1+\sqrt{\lambda_1}x_1,\ldots,y_{p_0}+\sqrt{\lambda_{p_0}}x_{2p_0-1})
  \quad ,
\end{eqnarray*}
where $X=\sum_{j=1}^p x_jX_j$.\\
It is easy to see that 
the representation
$\Pi$ is equivalent to $\rho_{r,\Lambda}$.
So $\Pi\in T_{rX^*+D_2(\Lambda)}$,
$\Pi$ is irreducible,
and its $K$-stability group denoted $K_\Pi$ 
equals $K_\rho$.\\
Let $\zeta_\alpha\in \Hc,\alpha\in E_l$ 
be given by:
$$
\zeta_\alpha:
\left\{
  \begin{array}{rcl}
    \Rb^{p_0}
    &\longrightarrow&
    \Rb\\
    y_1,\ldots,y_{p_0}
    &\longmapsto&
    \underset{j=1}{\overset{p_1}{\Pi}}
    h_{\alpha^j}(y_{m'_{j-1}+1},\ldots,y_{m'_j}) 
  \end{array}
\right. 
\quad,
$$
The vectors $\zeta_\alpha,\alpha\in E_l,l\in\Nb^{p_1}$ 
form an orthonormal basis of $\Hc$.\\
It can be proved that each vector space $V_l$ generated by
$\zeta_\alpha,\alpha\in E_l$
is $K_\Pi$-invariant.
Using relation between Laguerre and Hermite functions, 
we obtain:\\
\begin{lemma}
  \label{lem_fcnsph_repN}
  The spherical function associated 
  to $\Pi_{r,\Lambda}$ and $V_l$ is $\phi^{r,\Lambda,l}$.
\end{lemma}
\paragraph{Consequences for Sub-Laplacian.}
The (Kohn) sub-Laplacian is
$$
L
\,:=\,
- \sum_{i=1}^p X_i^2
\quad.
$$
It is a sub-elliptic $K$-invariant operator
(with analytic coefficients).
Consequently the spherical functions are eigenfunctions for this operator
(see Remark~\ref{rem_charac_sphfcn}).\\
Each representation $\Pi=\Pi_{r,\Lambda}$ 
induces the representation $d\Pi$
on the algebra of differential left invariant operators on $N_p$,
over the space of Schwartz functions $\Sc(\Rb^{p_0})$:
$$
\begin{array}{ccclr}
  j=1,\ldots,p_0
  &
  d\Pi(X_{2j-1})
  &=&
  \sqrt{\lambda_j} \partial_{y_j}
  &,
  \\
  j=1,\ldots,p_0
  &
  d\Pi(X_{2j})
  &=&
  i\sqrt{\lambda_j} y_j
  &,
  \\
  2p_0<j<p
  &
  d\Pi(X_j)
  &=&
  0
  &,
  \\
  j=p
  &
  d\Pi(X_p)
  &=&
  ir \Id
  &,
  \\
  \forall i<j 
  &
  d\Pi(X_{i,j})
  &=&
  \left\{
    \begin{array}{ll}
      i\lambda_{j'}\Id &\mbox{if}\;(i,j)=(2j'-1,2j') \\
      0&\mbox{otherwise} 
    \end{array}\right. 
  &. 
\end{array}
$$
For $L$, we get:
$$
d\Pi(L)
\,=\,
{r}^{2} \Id
-\sum_{i=1}^{p_0}
\lambda_i\left(\partial_{y_i}^2
  -y_i^2\right)
\quad.
$$
We recall
the Hermite function $ y=h_k$ satisfies the differential equation
$y''+(2k+1-x^2)y=0$ \cite[formula (5.5.2)]{szego}.
So, we get:
$$
d\Pi(L).\zeta_\alpha
\,=\,
\left( \sum_{j=1}^{p_1}
  \lambda_j (2 l_j+m_j)  
  +{r}^{2}\right)\zeta_\alpha
\quad,\quad \alpha\in E_l,\, l\in\Nb^{p_1}
\quad.
$$
We deduce (see Remark~\ref{rem_sphfcn_repN_radfcn}):
$$
L. \phi^{r,\Lambda,l}
\,=\,
\left( \sum_{j=1}^{p_1}
  \lambda_j (2 l_j+m_j) 
  +{r}^{2}\right)
\phi^{r,\Lambda,l}
\quad;
$$
this equality may also be computed directly using properties of the Laguerre functions.
\paragraph{Radial Plancherel measure.}
Here, we give the radial Plancherel measure.

Let $\Lc$ be  the set of   
$\Lambda=(\lambda_1,\ldots,\lambda_{p'})\in\Rb^{p'}$ such that
$\lambda_1\,>\, \ldots \,>\, \lambda_{p'}\,>\,0$.
We define the following measure on $\Lc$:
\begin{itemize}
\item $d\Lambda=d\lambda_1\ldots d\lambda_{p'}$ is the restricted Lebesgue measure on $\Lc$,
\item $\eta$ is the measure on $\Lc$ such that:
  $$
  d\eta(\Lambda)=
  \left\{\begin{array}{ll}
      c \Pi_{j<k} {(\lambda_j^2-\lambda_k^2)}^2 
      d\Lambda
      &\mbox{if}\; p=2p'\\
      c \Pi_i \lambda_i^2
      \Pi_{j<k} {(\lambda_j^2-\lambda_k^2)}^2 
      d\Lambda
      &\mbox{if}\; p=2p'+1
    \end{array}\right.
  \quad ,
  $$
  where the constant $c$ is chosen 
  in order to yield  the polar change of variables 
  over the space of antisymmetric matrices $\Ac_p$:
  $$    \int_{\Ac_p} g(A) dA \,=\,
  \int_{O_p} \int_{\Lc}
  g(k.D_2(\Lambda))
  d\eta (\Lambda) dk
  \quad.
  $$
\item $\eta'$ is
  the measure on $\Lc$ given by:
  $d\eta'(\Lambda)=\Pi_{i=1}^{p'} \lambda_id\eta(\Lambda)$,
\end{itemize}
Over $\Rb^+$, we define the measure $\tau$ given as
the Lebesgue measure if $p=2p'+1$, and the Dirac measure in 0
if $p=2p'$.

The (non radial) Plancherel measure 
is already known \cite[Section 6]{stric}.
With our notations,
it is  the measure $m$ given as the tensor product of
the Haar probability measure $dk$ on $K=O(p)$,
the measure $\eta'$ on $\Lc$,
and the measure $\tau$ on $\Rb^+$,
up to the constant $c(p)$  given by:
$$
c(p)
\,=\,
\left\{
  \begin{array}{ll}
    {(2\pi)}^{-\frac{p(p-1)}2+p'}
    &\quad
    \mbox{if}\; p=2p' \; ,\\
    2 {(2\pi)}^{-\frac{p(p-1)}2+p'-1}
    &\quad
    \mbox{if}\; p=2p'+1 \; .
  \end{array}\right.
$$
\begin{theorem}
  \label{thm_plancherel}
  $m$ is the (non radial) Plancherel measure;
  i.e. for $\psi\in L^2(N)$, we have:
  $$
  \nd{\psi}_{L^2(N)}^2
  \,=\,
  \int
  \nd{k.\Pi_{r,\Lambda}(\psi) }_{HS}^2
  dm(r,k,\Lambda)
  \quad.
  $$
  where $\nd{.}_{HS}$ denotes the Hilbert-Schmidt norm.
\end{theorem}
If we compute  the Hilbert-Schmidt square-norm 
of $k.\Pi_{r,\Lambda}(\psi)$
with the orthonormal basis $\{\zeta_\alpha,\alpha\in E_l, l\in\Nb^{p_1}\}$,
we deduce the radial Plancherel measure
(see Lemma~\ref{lem_fcnsph_repN}).
This is the measure, which we denote by $m^{\natural}$,
given as  the tensor product
of $\eta'$ on $\Lc$, 
and the counting measure $\sum$ on $\Nb^{p'}$,
and  the measure $\tau$ on $\Rb^+$,
up to the normalizing constant $c(p)$:
\begin{theorem}
  \label{thm_mplanch}
  $m^{\natural}$ is the radial Plancherel measure for $(N_p,O_p)$,
  i.e. for a $K$-invariant function $\psi\in L^2(N)$, we have:
  $$
  \nd{\psi}_{L^2(N)}^2
  \,=\,
  \int
  \nn{<\psi,\phi^{r,\Lambda,l}>}^2
  dm^{\natural}(r,\Lambda,l)
  \quad.
  $$
\end{theorem}
We can also compute directly the radial Plancherel measure 
$m^\natural$, using the properties 
of Laguerre functions and Euclidean Fourier transform
(see \cite{moi}).

\end{document}